\documentclass[12pt]{article}
\usepackage{amsmath}
\usepackage{graphicx,psfrag,epsf}
\usepackage{enumerate}
\usepackage{natbib}
\usepackage{amsmath,amsthm,amsfonts,epsfig,amssymb,natbib,eucal,eufrak}
\usepackage{booktabs}
\usepackage{multirow}
\usepackage{siunitx}
\usepackage{qtree}
\usepackage{hyperref}

\usepackage{float}
\restylefloat{table}
\usepackage{color}

\newtheorem{cor}{Corollary}

\newcommand{\blind}{0}

\begin{document}

\def\spacingset#1{\renewcommand{\baselinestretch}%
{#1}\small\normalsize} \spacingset{1}


\if0\blind
{
  \title{\bf A note on the closed-form solution for the longest head run problem of Abraham de Moivre}
  \author{Yaakov Malinovsky
    \thanks{email: yaakovm@umbc.edu}
   \\
    Department of Mathematics and Statistics\\ University of Maryland, Baltimore County, Baltimore, MD 21250, USA\\
}
  \maketitle
} \fi

\if1\blind
{
  \bigskip
  \bigskip
  \bigskip
  \begin{center}
    {\LARGE\bf Title}
\end{center}
  \medskip
} \fi

\bigskip
\begin{abstract}
The problem of the longest head run was introduced and solved by Abraham de Moivre in the second edition of his book {\it Doctrine of Chances} \citep{D1756}. The closed-form solution as a finite sum involving binomial coefficients was provided in \cite{U1937}.
Since then, the problem and its variations and extensions have found broad interest and diverse applications. Surprisingly, a very simple closed form can be obtained, which we present in this note.

\end{abstract}

\noindent%

\spacingset{1.45} 
\section{Introduction}
In a series of $n$ independent trials, an event $E$ has a probability $p$ of occurrence for each trial.
If, in these trials, event $E$ occurs {\it at least} $r$ times without interruption, then we have a {\it run} of size $r$.
What is the probability $y_{n,r}$ of having a run of size $r$ in $n$ trials?
This problem was formulated and solved by Abraham de Moivre in the second edition of his book {\it The Doctrine of Chances: or, A Method of Calculating the Probabilities of Events in Play} (\cite{D1756}, Problem LXXXVIII, p. 243).
Although more than 280 years have passed since then, de Moivre's problem and its variations remain of great interest in probability and statistics; see for example \cite{N2017} and references therein.

Formally, let $X_1,\ldots,X_n$ be a sequence of independent and identically distributed Bernoulli random variables with the probability of success $p=P(X_1=1)=P(E)$ and the probability of failure $1-p=P(X_1=0)$.
We denote by
$$L_{n}=\max\left\{k: X_{i+1}=\cdots=X_{i+k}=1,\,i=0,1,\ldots,n-k\right\}$$
the length of the longest consecutive series of trials where event $E$ occurs.
Then, recalling that having a run of size $r$ is defined as the occurrence of at least $r$ consecutive successes
in $n$ iid Bernoulli trials with corresponding probability $y_{r,n}$, we obtain

\begin{equation*}
\label{eq:d}
y_{n,r}=P\left(L_{n}\geq r\right).
\end{equation*}

de Moivre did not provide a proof, but demonstrated a method of finding $y_{n,r}$. Reviewing that method, one can see that he used the method of generating functions. He demonstrated the method with an example of ten trials having $p=1/2$, in which
the probability of a run of size $3$ equals $65/128$.

The closed-form solution of $y_{n,r}$ was given by \cite{U1937} as a polynomial with binomial coefficients, arrived at by
first obtaining a difference equation and then using a method of generating functions to solve it.
Surprisingly, a simple closed form of $y_{n,r}$ can be obtained as a corollary from the difference equation given by Uspensky.
This closed-form solution seems to have never been reported in the literature. In this note, we present it along with Uspensky's original derivations.

\section{Uspensky's solution}
We present Uspensky's solution (\cite{U1937}, pages 77-79) while keeping his original notations.
This solution demonstrates the power of the use of ordinary linear difference equations along with the generating functions.
Let's consider $n+1$ trials with the corresponding probability $y_{n+1,r}$.
A run of size $r$ in $n+1$ trials can happen in two mutually-exclusive ways:
\begin{itemize}
\item[(W1)]: the run is obtained in the first $n$ trials or
\item[(W2)]: the run is obtained as of trial $n+1$.
\end{itemize}
(W2) means that among the first $n-r$ trials, there is no run of size $r$; event $E^{C}$ occurred at trial $n-r+1$;
and event $E$ occurred in the trials $n-r+2,\ldots,n+1$.
Combining (W1) and (W2), we obtain a linear difference equation of order $r+1$,
\begin{equation}
\label{eq:R}
y_{n+1,r}=y_{n,r}+(1-y_{n-r,r})q p^r,
\end{equation}
with initial conditions
\begin{equation*}
y_{0,r}=y_{1,r}=\cdots=y_{r-1,r}=0,\,\,\,y_{r,r}=p^r,
\end{equation*}
where $q=1-p$.

Using a method of generating functions, Uspensky obtained the closed-form solution of \eqref{eq:R} 
\begin{align}
\label{re:U}
&
y_{n,r}=1-\beta_{n,r}+p^{r}\beta_{n-r,r}\nonumber \\
&
\beta_{n,r}=\sum_{l=0}^{\left \lfloor{\frac{n}{r+1}}\right \rfloor }(-1)^{l}{n-lr \choose l}(qp^{r})^{l},
\end{align}
where $\lfloor{x} \rfloor$ is defined as the greatest integer less than or equal to $x$. 

Going back to de Moivre's original example, where $n=10,r=3$ and $p=1/2$, and using \eqref{re:U}, we obtain $y_{n,r}=65/128$, which coincides with de Moivre's answer on page 245 of his {\it Doctrine of Chances} \citep{D1756}.

\section{Closed-form solution for the case $r\geq n/2$}
Surprisingly, a simple closed-form solution follows from Uspensky's equation \eqref{eq:R} as follows:
Replacing $n+1$ with $n$, the difference equation \eqref{eq:R} becomes
\begin{equation}
\label{eq:Y}
y_{n,r}=y_{n-1,r}+(1-y_{n-1-r,r})q p^r,
\end{equation}
with initial conditions $y_{0,r}=\cdots=y_{r-1,r}=0, y_{r,r}=p^r$.

If $n-1-r\leq r-1$ (i.e. $n/2\leq r$), then $y_{n-1-r,r}$=0 and we obtain from \eqref{eq:Y}
\begin{equation}
\label{eq:YY}
y_{n,r}=y_{n-1,r}+q p^r\,\,\,\,\text{for}\,\,\,\,\,n/2\leq r,
\end{equation}
with initial conditions  $y_{0,r}=\cdots=y_{r-1,r}=0, y_{r,r}=p^r$.

Iterating \eqref{eq:YY}, we obtain the following corollary:

\begin{cor}
\label{cor:U}
If $n/2\leq r \leq n$, where $r$ is an integer, then
\begin{equation}
\label{eq:new}
y_{n,r}=p^{r}+(n-r)p^rq.
\end{equation}
\end{cor}

An alternative explanation for \eqref{eq:new} is that, when $n/2\leq r \leq n$, we cannot have two separate runs of size $r$. Therefore, a run of size $r$ in $n$ trials either begins at position 1 (this happens with probability $p^{r}$) or at the positions $2,\ldots,n-(r-1)$ (at each such position, this happens with probability $qp^r$). Then, denote by $R_i$ an event in which a run of size $r$ begins at position $i, i=1,2\ldots,n-(r-1)$. Since the events $R_1,R_2,\ldots, R_{n-(r-1)}$ are disjoint, we obtain that $y_{n,r}=P(R_1)+P(R_2)+\ldots+P(R_{n-(r-1)})=p^{r}+(n-r)qp^r$.

\section{Comments}
There are a number of interesting problems discussed in the \cite{U1937} book,
many of which have roots in the classics of probability, their origins tracing back to founders of modern-probability such as Pascal, Fermat, Huygens, Bernoulli, de Moivre, Laplace, Markoff, Bernstein, and others.
For example, \cite{U1937} considered another problem of de Moivre's that was later discussed and extended by \cite{DZ1991}. A large collection of classic problems in probability with historical comments and original citations are nicely presented in the book by \cite{G2012}.

There are many follow-ups on and extensions of de Moivre's longest head run problem.
An interesting recursive solution of the problem in the case of the fair coin was given by \cite{ST1979}.
That problem was then extended to the Markov chain setting, where Uspensky's generation function was generalized to the case of dependent observations \citep{N1989}.
The problem has also found applications in numerous fields. Among which are reliability \citep{DLR1982}, computational biology \citep{S2000}, and finance where time dependent-sequences naturally occur (see \cite{N2011} and references therein).

\section*{Acknowledgement}
I thank Serguei Novak for comments on the early version of this article, and
I thank Danny Segev for following the early version of this article on arXiv and suggesting the alternative proof that is now presented after Corollary 1.
I also would like to thank the referee for insightful and helpful comments that led to significant improvements in the paper.

{}

\end{document}